\documentclass[a4paper,11pt]{article}
\usepackage{cite}
\usepackage{graphicx}
\usepackage[tight,footnotesize]{subfigure}
\usepackage{amsfonts,latexsym,amssymb,amsmath}
\usepackage{epstopdf}
\usepackage{tikz}
\usetikzlibrary{arrows}
\newtheorem{theorem}{Theorem}[section]

\newtheorem{proposition}[theorem]{Proposition}
\newtheorem{algorithm}[theorem]{Algorithm}
\newtheorem{remark}{Remark}
\newtheorem{definition}{Definition}
\begin{document}

\title{Finding the Exact Delay Bound for Consensus of Linear Multi-Agent Systems}
\author{Rudy Cepeda-Gomez\thanks{Institute of Automation, University of Rostock, Rostock, Germany. Email: rudy.cepeda-gomez@uni-rostock.de}}

\maketitle

\begin{abstract}
This paper focuses on consensus problems for high-order, linear multi-agent systems. Undirected communication topologies and fixed, uniform communication time delay are taken into account. This class of problems has been widely study in the literature, but there are still gaps concerning the exact delay stability bounds in the domain of the delays. The more common analysis employed is based on Lyapunov-Krasowskii functionals, which produce very conservative results that are cumbersome to apply. As an alternative, we employ the Cluster Treatment of Characteristic Roots paradigm to study the stability of the system in the space of the delay. This allows the generation of exact and exhaustive delay bounds in an efficient manner. Before the stability analysis, a state transformation is performed to decouple the system and simplify the problem, as it was previously done for consensus problem of agents with simpler dynamics. Simulation results are presented to support the analytical claims. 
\end{abstract}

\section{Introduction}
In the past two decades considerable research effort has been placed in the multi-agent systems area. Within this field, the particular topic of consensus has received a great deal of attention. In a consensus problem the member of a team aim at reaching an agreement among them. This problem is treated by \cite{Vicsek1995} in the early on, where agents try to align their headings using discrete-time representation. Later, \cite{Olfati-Saber2004} introduced a formal framework for the consensus problem for agents driven by first-order dynamics. That study considers directed communication topologies, both for fixed and switching types, and introduces some graph theoretical results useful for the stability analysis of such agreement protocols. Under the simplifying features of first-order dynamics, they also consider time-delayed communications in the case of a fixed topology. Many other researchers suggest further extensions to this work, proposing consensus protocols for agents driven by second-order dynamics \cite{Ren2008,Ren2005} without time delays in the communication channels.

The works by Lin and his co-workers include the study of agents driven by second-order discrete \cite{Lin2009-2} and continuous \cite{Lin2009-1} dynamics, and applications of different control laws, ranging from a basic P-D-like (proportional and Derivative) control logic \cite{Lin2008-1} to $H_\infty$ structure \cite{Lin2008-2}. Most of these works entail time-delayed communications and directed topologies. Sun and Wang also consider the consensus problem with fixed \cite{Sun2009-2} and time-varying delays \cite{Sun2009-1}. \cite{Peng2009} study a leader follower case, also with time-varying delays. \cite{Meng2011} present a comparison study among several different swarm control laws, including first- and second-order agents and considering two time delays. A more detailed review of the recent results and open challenges in this topic can be found in the work of \cite{Cao2013b}.

Most of the works mentioned study systems with very simple dynamics. On the contrary, \cite{Zhang2014} is one of the few works that address the problem of state consensus on networks of agents driven by higher-order, general linear dynamics under undirected communications and affected by constant time delays. By transforming the consensus problem into a robust stability problem for an uncertain system, they arrive to an allowable delay bound for systems with fixed and switching interconnection topologies.

When it comes to the stability analysis of consensus systems with time delays, almost all the previous works rely on methodologies based on Lyapunov/Krasovskii or Razhumikin theorems  \cite{Lin2008-1,Lin2008-2,Lin2009-1,Lin2009-2,Meng2011,Zhang2014} or on the generalised Nyquist criterion \cite{Liu2009,Muenz2010}. All of these treatments provide only sufficient conditions on the delays to achieve asymptotic stability. They produce only very conservative results, leading to stability bounds for very small delays. In addition, because these results are based on the solutions of some Linear Matrix Inequalities (LMI), they are always imprecise, conservative and, very critically, cumbersome to deploy.

As a different way to analyse the consensus problem with respect to the delay, \cite{TAC,IJC2} introduced a methodology for the analysis of consensus protocols with single and multiple, rationally independent time delays. This earlier work is based on the combination of a simplifying factorization procedure over the characteristic equation of the system and the deployment of a crucial stability paradigm, which is called the Cluster Treatment of Characteristic roots CTCR \cite{Ergenc2007,Fazelinia2007,Olgac2005}. The CTCR provides a tool for assessing the stability of linear time-invariant (LTI) systems with multiple rationally independent delays. This new method uniquely creates exact, exhaustive and explicit stability regions in the domain of the delays. 

All of our earlier investigations \cite{TAC,IJC1,IJC2} have focused on agents driven by very simple, double integrator dynamics. In the present work, we deploy our methodology of state decomposition and stability analysis using CTCR to the problem of state consensus for higher order, linear systems. We show that this methodology provides \emph{exact} stability bounds in the domain of the delay, and its deployment is straight forward. We compare our results to those of \cite{Zhang2014} to show the advantages of our method.

The rest of the paper is organized as follows. Section \ref{sec:problem} presents the system under study and the state transformation that simplifies its stability analysis. Section \ref{sec:ctcr} presents the stability analysis with respect to the time delay, including a revision of the CTCR paradigm. Simulation results to support the analytical claims are presented in section \ref{sec:simulations}, whereas Section \ref{sec:switch} makes some comments on the switching topologies case. Finally, concluding remarks and directions for future research are given in section \ref{sec:conclusions}. There is also an appendix that introduces the Kronecker summation operation and some of its properties.
\section{Problem Statement}\label{sec:problem}
In this paper, we consider a network of $n$ agents, each one with continuous-time linear dynamics given by
\begin{equation}
\dot{\mathbf{x}}_i\left(t\right)=\mathbf{A}\mathbf{x}_i\left(t\right)+\mathbf{B}\mathbf{u}_i\left(t\right),
\label{eq:agentdyn}
\end{equation} 
where $\mathbf{x}\in\mathbb{R}^{p}$ being the state of the agent and $\mathbf{u}\in\mathbb{R}^{q}$ its control input. The matrices $\mathbf{A}\in\mathbb{R}^{p\times p}$ and $\mathbf{B}\in\mathbb{R}^{p\times q}$ are assumed to be constant. Furthermore, $\mathbf{A}$ is assumed to be not Hurwitz, to prevent the agents from reaching consensus at $\mathbf{x}_i=0$ without cooperation. We state that the agents reach consensus if $\lim_{t\rightarrow\infty}\left\|\mathbf{x}_i-\mathbf{x}_k\right\|=0$ for any $i,\,k\,\in\,[1,n]$.

In order to reach consensus, the $n$ agents exchange state information through bidirectional communication channels. Their communication topology, therefore, can be described by an undirected graph. The adjacency matrix of this graph is denoted by $\mathcal{A}=\left[a_{ik}\right]\in\mathbb{R}^{n\times n}$. It has elements $a_{ik}=1$ if there is a communication channel among agents $i$ and $k$ and $a_{ik}=0$ otherwise. The diagonal elements are assumed as zero, i.e., $a_{ii}=0$, for $i=1,\,2,\ldots,\,n$. The Laplacian matrix of the communication topology, $\mathbf{L}=\left[l_{ik}\right]\in\mathbf{R}^{n\times n}$, is defined such that its diagonal elements are equal to sum of the corresponding row of the adjacency matrix, $l_{ii}=\sum_{k=1}^{n}{a_{ik}}$, and its off-diagonal elements are the negative of the elements in $\mathcal{A}$, $l_{ik}=-a_{ik}$.

It is assumed that in every communication channel there is a time delay, $\tau$, which is constant and uniform through the network. With this delay, the agents implement a control law given by:
\begin{equation}
\mathbf{u}_i\left(t\right)=\mathbf{K}\sum_{k=1}^{n}{a_{ik}\left(\mathbf{x}_k\left(t-\tau\right)-\mathbf{x}_i\left(t-\tau\right)\right)},
\label{eq:control}
\end{equation} 
in which $\mathbf{K}\in\mathbb{R}^{q\times p}$ is a constant gain matrix. 

Using this control law, the dynamics of the complete system of $n$ agents is given by:
\begin{equation}
\dot{\mathbf{x}}\left(t\right)=\left(\mathbf{I}_n\otimes\mathbf{A}\right)\mathbf{x}\left(t\right)-\left(\mathbf{L}\otimes\mathbf{BK}\right)\mathbf{x}\left(t-\tau\right),
\label{eq:fullsystem}
\end{equation}
where the full state vector is formed by the concatenation of the states vectors of the individual systems, $\mathbf{x}\left(t\right)=\left[\mathbf{x}_1^{\rm{T}}\ \mathbf{x}_2^{\rm{T}}\ \cdots\ \mathbf{x}_n^{\rm{T}}\right]^{\rm{T}}$, $\mathbf{I}_n$ is the identity matrix of size $n\times n$ and $\otimes$ represents the Kronecker product \cite{kron}. 

A well established practice in the consensus literature is the transformation of the full dynamics \eqref{eq:fullsystem} into a stability analysis problem. A straightforward approach that decomposes \eqref{eq:fullsystem} into a set of reduced order system with identical structure but different parametric values, was introduced by \cite{TAC} and is used in this work. This method is based on the diagonalization of the Laplacian matrix.

Since the Laplacian is a symmetric matrix, there is an orthogonal matrix $\mathbf{T}\in\mathbb{R}^{n\times n}$ such that $\boldsymbol{\Lambda}=\mathbf{T}^{\rm{T}}\mathbf{L}\mathbf{T}$. Here $\boldsymbol{\Lambda}\in\mathbb{R}^{n\times n}$ is a diagonal matrix that contains the eigenvalues of $\mathbf{L}$, which are all real. A known fact from algebraic graph theory \cite{Godsil} states that the smallest eigenvalue of $\mathbf{L}$ is zero whenever the communication topology is connected. According to this and without loss of generality, we sort the eigenvalues of the Laplacian as $0=\lambda_1<\lambda_2\leq\lambda_3\leq\cdots\leq\lambda_n$. The $i$-th column of the matrix $\mathbf{T}$, denoted as $\mathbf{t}_i\in\mathbb{R}^{n}$ is the normalized eigenvector corresponding to the $i$-th eigenvalue. Following the theory, this implies that when the topology is connected, $\mathbf{t}_1$, the first column of $\mathbf{T}$, equals $1/\sqrt{n}\mathbf{1}_n$, where $\mathbf{1}_n$ is a $n\times 1$ vector with all its  elements equal to one.

Using the definitions and results of the previous paragraph, a state transformation is defined using $\mathbf{T}$:
\begin{equation}
\boldsymbol{\xi}=\left(\mathbf{T}^{\rm{T}}\otimes\mathbf{I}_q\right)\mathbf{x}.
\label{eq:trans}
\end{equation}
When \eqref{eq:trans} is used in \eqref{eq:fullsystem}, after a few manipulations using the properties of the Kronecker product, we arrive to:
\begin{equation}
\dot{\boldsymbol{\xi}}\left(t\right)=\left(\mathbf{I}_n\otimes\mathbf{A}\right)\boldsymbol{\xi}\left(t\right)-\left(\boldsymbol{\Lambda}\otimes\mathbf{BK}\right)\boldsymbol{\xi}\left(t-\tau\right).
\label{eq:diagsystem}
\end{equation}
Since $\boldsymbol{\Lambda}$ is a diagonal matrix, the system in \eqref{eq:diagsystem} is actually composed by $n$ decoupled subsystems of order $q$ with dynamics of the form:
\begin{equation}
\dot{\boldsymbol{\xi}}_i\left(t\right)=\mathbf{A}\boldsymbol{\xi}_i\left(t\right)-\lambda_i\mathbf{BK}\boldsymbol{\xi}_i\left(t-\tau\right),
\label{eq:subsys}
\end{equation}
where the state of the individual subsystems is $\boldsymbol{\xi}_i=\left(\mathbf{t}_i^{\rm{T}}\otimes\mathbf{I}_q\right)\mathbf{x}$. 

Paying close attention to the properties of the Laplacian and the definition of the subsystems, we can see that for $i=1$ we have $\boldsymbol{\xi}_1=1/\sqrt{n}\sum_{i=1}^n\mathbf{x}_i$ and:
\begin{equation}
\dot{\boldsymbol{\xi}}_1\left(t\right)=\mathbf{A}\boldsymbol{\xi}_1\left(t\right).
\label{eq:centroid}
\end{equation}
This leads to the conclusions that the first subsystem dictates the dynamics of the average of the state of the members of the group and that it is not affected by the time delay. Indeed, this subsystem dictates the group decision value, i.e., the value of the states upon which the agents agree if they reach consensus.

The remaining $n-1$ subsystems in \eqref{eq:subsys} are related to the dynamics of different linear combinations of the states of the original agents. They represent the so called \emph{disagreement dynamics}. These systems must be all simultaneously stable for the agents to reach consensus. These facts have been formally stated and proven for several different consensus protocols. The interested reader can look at the works of \cite{IJC1}, \cite{Olfati-Saber2004}, \cite{Sun2009-1}, \cite{Zhang2014} and others.

In the work of \cite{Zhang2014} an uncertain system is defined based on the distribution of the $n-1$ non-zero eigenvalues of $\mathbf{L}$. Using a Lyapunov-Krasowskii type of analysis, they find a delay stability bound that guaranties stability for the uncertain system. This, in turns, guarantees stability for all the $n-1$ disagreement subsystems in \eqref{eq:subsys}. This approach, however, has several disadvantages. First, the delay bound found is extremely conservative. This conservatism is introduced by the consideration of the uncertain system and by the very nature of the Lyapunov-Krasowskii approach. The second disadvantage is that the delay bound is not easy to calculate. In order to find the bound, several matrices must be defined, and they mus satisfy a set of Linear Matrix Inequalities. Finding such matrices is not an easy task. The third main disadvantage is connected to this. For the same system, changing the definition of one of the matrices involved in the computation of the delay bound can change its value. Therefore, the delay bound is not unique for a given system.

As a way to overcome these difficulties, in the following section we present an analysis methodology based on the Cluster Treatment of Characteristic Roots paradigm. This is a very efficient approach that allows us to find, in an explicit manner, the \emph{exact} stability bound for a given system in a very short time.

\section{Stability Analysis Using CTCR}\label{sec:ctcr}
\subsection{The CTCR Paradigm}
Consider a linear, time-invariant, time delay system (\textsc{lti-tds}) with dynamics:
\begin{equation}
\dot{\mathbf{x}}\left(t\right)=\mathbf{Ax}\left(t\right)+\mathbf{Bx}\left(t-\tau\right).
\label{eq:ltitds}
\end{equation} 
It is well known that the characteristic equation of this class of systems, given by:
\begin{equation}
{\rm det}\left(s\mathbf{I}-\mathbf{A}-\mathbf{B}{\rm e}^{-s\,\tau}\right)=0,
\label{eq:ce}
\end{equation}
exhibits an infinite number of roots, due to the transcendental term ${\rm e}^{-s\,\tau}$ introduced by the delay. For the system \eqref{eq:ltitds} to be stable, all those roots must have negative real part. To assure this, is a complex task.

The Cluster Treatment of Characteristic Roots, CTCR, paradigm \cite{Olgac2002,Olgac2005,Olgac2006} is a mathematical tool that enables us to solve this task in an exact and explicit manner. It builds upon the \emph{D-subdivision method}, also known as the \emph{continuity argument} \cite{kolmanovski1986}. This argument states that for a retarded system (a system in which the highest derivative of the state is not affected by the delay) like \eqref{eq:ltitds} the solutions of \eqref{eq:ce} are a continuous function of the delay. Because of this, a change in the stability posture of the system with respect to the delay is only possible when, for a given value of the delay, $\tau_c$, the system exhibits a conjugated pair of purely imaginary characteristic roots, $s=\pm\rm{j}\omega_c$. Therefore, the delay space, $\tau\in\mathbb{R}^+$ is divided in \emph{pockets} in which the number of unstable characteristic roots of the system remains constant.

Another important and well known fact is that the purely imaginary roots of \eqref{eq:ce} are periodic with respect to the delay. This means that is $s=\rm{j}\omega_c$ is a root of \eqref{eq:ce} for $\tau=\tau_c$, the same root appears when $\tau=\tau_c+2\pi/\omega_c$.

With this two facts at hand, we present three important definitions, before moving to the main propositions of CTCR.

\begin{definition}
\label{def:kernel}
\emph{Kernel} $\wp_0$: The points on the positive real line $\mathbb{R}^{+}$ that consist exhaustively of all the points $\tau$ values which cause a pair of imaginary characteristic roots of \eqref{eq:ce} at $s=\pm\rm{j}\omega_c$ and satisfy the constraint $0<\tau\omega_c<2\pi$, are called the \emph{kernel points}. These are the smallest possible delay values that create the given pair of imaginary roots at the frequency $\omega_c$.
\end{definition}
\begin{definition}
\label{def:offspring}
\emph{Offspring} $\wp$: The points obtained from the kernel by using the periodicity of the imaginary roots with respect to the delay: $\tau=\tau+2k\pi/\omega_c$ are called the offspring.
\end{definition}
\begin{definition}
\label{def:rt}
\emph{Root Tendency}, $RT$: At any point $\tau\in\wp_0\cup\wp$ an infinitesimal increase in the delay creates a transition of the root. Such transition can be to the right or to the left half of the complex plane. The Root Tendency, $RT$, indicates the direction of this transition:
\begin{equation}
\left.RT\right|_{s=\rm{j}\omega}=\mathbf{sgn}\left[\Re\left(\left.\frac{\partial s}{\partial\tau}\right|_{s=\rm{j}\omega}\right)\right]
\label{eq:rt}
\end{equation}
Clearly these root tendencies are $-1$ for stabilizing and $+1$ for destabilizing root crossings on the imaginary axis.
\end{definition}

The two key propositions of CTCR follow.

\begin{proposition}
\label{p:prop1}
\emph{Finite Number of Kernel Points} A given LTI-MTDS can exhibit only a finite number, $m$, of kernel points. This number is upperbounded by the square of the order of the system: $m<n^2$ \cite{Ergenc2007}. 
\end{proposition}
\begin{proposition}
\label{p:prop2}
\emph{Invariance of the Root Tendency} Take an imaginary characteristic root, $s=\rm{j}\omega_c$, caused by any one of the infinitely many grid points in the kernel and offspring sets. The root tendency of these imaginary roots remains invariant so long as offspring points are obtained from the same point in the kernel. That is, the root tendency with respect to the variations of $\tau$ is invariant from the kernel to the corresponding offspring, $\tau=\tau+2k\pi/\omega_c$. 
\end{proposition}
\begin{remark}
It is important to highlight that the definitions and propositions presented here are tailored for a system with single time delay, as \eqref{eq:ltitds}. CTCR, however, is a very powerful tool that can handle linear systems with any number of delays (within computational feasibility). The full definitions and proofs for these cases can be found in the works of \cite{Olgac2005}, \cite{Fazelinia2007} and \cite{Ergenc2007}. 
\end{remark}

Jointly the kernel and the offspring sets contain all the points in the delay space where the system in \eqref{eq:ce} has imaginary characteristic roots and therefore can change its stability posture. This complete set constitutes the departure point of the CTCR paradigm. Once this set is known, the number of unstable roots of the system at $\tau=0$ is found. This is the number of unstable roots for the first pocket. Then, moving along the real line, the number of unstable roots is increased by two when a point in $\wp_0\cup\wp$ with $RT=1$ is reached or decreased by two when such a point has $RT=-1$. The pockets in which the number of unstable roots is zero are declared as stable operating regions for the system.

It is crucial that an appropriate computational tool is used to capture all the kernel points exhaustively. Several different methods have been developed for this task. A simple algebraic approach can be used when the system has a single time delay with no commensurate (integer multiples of the delay) elements, as it was done for the Delay Resonator active vibration absorber \cite{Olgac1994} as well as in previous applications of CTCR to consensus problems with single delay \cite{IJC1,TAC,EJC}. The Rekasius substitution can be used when there are commensurate delays \cite{Olgac2002,Olgac2006} or even multiple delays \cite{Olgac2005}. The \emph{Building Block} concept \cite{Fazelinia2007} is also useful to simplify the task when there are multiple delays, and it can be deployed using the Rekasius substitution \cite{Fazelinia2007} or a half-angle tangent substitution, as used in previous works concerning consensus under multiple time delays \cite{IJC2,SCL}. 

Because of its easy implementation, the Kronecker summation method introduced by \cite{Ergenc2007} is used in this study to find the stability regions of the subsystems \eqref{eq:subsys} for $i=2,\,3\ldots,\,n$. This analysis is presented in the following subsection. 
\subsection{Analysis of the Individual Systems}
The systems in \eqref{eq:subsys} have the generic form:
\begin{equation}
\dot{\mathbf{x}}\left(t\right)=\mathbf{A}\mathbf{x}\left(t\right)+\mathbf{C}\mathbf{x}\left(t-\tau\right),
\label{eq:gensys}
\end{equation}
with $\mathbf{C}=-\lambda_i\mathbf{BK}$. The characteristic equation of this system has the form:
\begin{equation}
{\rm det}\left(s\mathbf{I}-\mathbf{A}-\mathbf{C}z\right)=0,
\label{eq:ce2}
\end{equation}
where $z={\rm e}^{-s\,\tau}$. Finding the roots of this equation is equivalent to finding the eigenvalues of the matrix:
\begin{equation}
\mathbf{D}_1=\mathbf{A}+\mathbf{C}z.
\label{eq:eigenvalue1}
\end{equation}
Given that the system matrices $\mathbf{A}$ and $\mathbf{C}$ are real valued, the solutions to \eqref{eq:ce2} come in complex conjugate terms. This implies that \eqref{eq:ce2} is also satisfied by the complex conjugates of $s$ and $z$:
\begin{equation}
{\rm det}\left(s^{\ast}\mathbf{I}-\mathbf{A}-\mathbf{C}z^{\ast}\right)=0.
\label{eq:ceconj}
\end{equation}  
This equation, in turn, is equivalent to determining the eigenvalues of:
\begin{equation}
\mathbf{D}_2=\mathbf{A}+\mathbf{C}z^{\ast}.
\label{eq:eigenvalue2}
\end{equation}

Since we are interested in finding the kernel and offspring sets, i.e. the purely imaginary characteristic roots of the system, we focus on those solutions to \eqref{eq:ce2} that have $s={\rm j}\omega$, and $z={\rm e}^{-{\rm j}\omega\,\tau}$. Then $z$ is a unitary complex number and $z^{\ast}=z^{-1}$.

When there is a common solution to \eqref{eq:ce2} and \eqref{eq:ceconj} $s={\rm j}\omega$ is an eigenvalue of $\mathbf{D}_1$ while $\mathbf{D}_2$ has $s^{*}=-{\rm j}\omega$ as an eigenvalue. Obviously, the sum of these two eigenvalues is zero, which implies that the Kronecker summation $\mathbf{D}_1\oplus\mathbf{D}_2$ must have one eigenvalue equal to zero (see appendix \ref{sec:appendix} and the works of \cite{Brewer1978} and \cite{Bernstein2005} for a definition and properties of the Kronecker summation operation). Since $\mathbf{D}_1\oplus\mathbf{D}_2$ has a zero eigenvalue, its determinant must be zero, and this leads to the definition of the \emph{auxiliary characteristic equation} for system \eqref{eq:gensys}:
\begin{equation}
ACE(z)={\rm det}\left[\mathbf{D}_1\oplus\mathbf{D}_2\right]={\rm det}\left[\left(\mathbf{A}+\mathbf{C}z\right)\oplus\left(\mathbf{A}+\mathbf{C}z^{-1}\right)\right]=0.
\label{eq:ace}
\end{equation}
Notice that \eqref{eq:ace} is an algebraic equation in $z$, which can be efficiently solved by several numerical tools. The order of this equation depends on the rank of $\mathbf{C}$, but can be at most ${\rm rank}\left(\mathbf{C}\right)^2$. This is a proof of proposition \ref{p:prop1}.

From the solutions $z$ to \eqref{eq:ace} we extract those that have unit norm. This results are then used in \eqref{eq:ce2} to find the values of $\omega$, and from $z={\rm e}^{-{\rm j}\omega\tau}$ the values of the delay are obtained. When finding the delays, one must be sure to obtain the minimum positive value, to make sure it corresponds to a kernel point. The root tendencies are obtained from they definition in \eqref{eq:rt} and this completes the construction of the kernel. The periodicity of the roots and the D-subdivision concept are then applied to find the stable and unstable operating regions in the space of the delay.

The stability analysis method is described in the following algorithm.
\begin{algorithm}
For a system in the form \eqref{eq:gensys}.
\begin{enumerate}
\item Define the auxiliary variable $z$.
\item Construct the auxiliary characteristic equation, $ACE(z)$, in \eqref{eq:ace}, and find its solutions.
\item From the solutions to \eqref{eq:ace} extract those with $\left|z\right|=1$.
\item Replace these values of $z$, together with $s={\rm j}\omega$ in \eqref{eq:ce2} and solve for $\omega$.
\item Using $z={\rm e}^{-{\rm j}\omega\tau}$ and the periodicity property, find the kernel delay points.
\item Using \eqref{eq:rt} find the root tendency.
\item Extend the kernels to obtain the offspring.
\item Start from the roots of the system for $\tau=0$ and deploy the D-Sub-division concept to obtain the full stability map in the space of the delay.
\end{enumerate}
\end{algorithm}

The repeated application of the previous algorithm to the $n-1$ individual disagreement systems of the form \eqref{eq:subsys} can be performed very efficiently. Once the stability regions of each subsystem are found, they are intersected to obtain the stability regions for the complete system, and the \emph{exact} delay bound of the consensus system.

\section{Numerical Examples}\label{sec:simulations}
\begin{figure}
\centering
\begin{tikzpicture}
\node (n1) [draw,circle] at( 0, 1.5){1};
\node (n2) [draw,circle] at( 1.5, 0){2};
\node (n3) [draw,circle] at( 1.5,-1.5){3};
\node (n4) [draw,circle] at(-1.5,-1.5){4};
\node (n5) [draw,circle] at(-1.5,0) {5};
\draw [thick,stealth-stealth] (n1.east)  to [out=0,in=90]   (n2.north);
\draw [thick,stealth-stealth] (n2.south) to                 (n3.north);
\draw [thick,stealth-stealth] (n3.west)  to                 (n4.east);
\draw [thick,stealth-stealth] (n4.north) to                 (n5.south);
\draw [thick,stealth-stealth] (n5.north) to [out=90,in=180] (n1.west);
\end{tikzpicture}
\caption{Ring topology with five agents.}
\label{fig:ringtopology}
\end{figure}
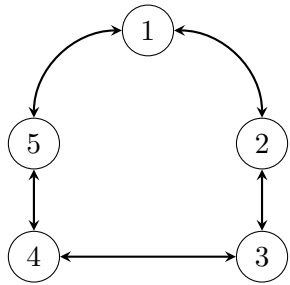

We test our methodology on the same example presented by \cite{Zhang2014}. Consider 5 agents driven by \eqref{eq:agentdyn} with:
\begin{displaymath}
\mathbf{A}=\left[\begin{array}{crc}0.2&0&0\\0&0&1\\1&-1&0\end{array}\right]\quad
\mathbf{B}=\left[\begin{array}{cc}1&0\\0&1\\1&0\end{array}\right],
\end{displaymath}
interacting under the communication topology of Figure \ref{fig:ringtopology}.
The adjacency matrix of this topology is given by:
\begin{equation}
\mathcal{A}=\left[\begin{array}{ccccc}0&1&0&0&1\\1&0&1&0&0\\0&1&0&1&0\\0&0&1&0&1\\1&0&0&1&0\end{array}\right],
\label{eq:ag}
\end{equation}
and its corresponding Laplacian matrix has three different eigenvalues: $0$, $1.3820$ and $3.6180$, the last two of them have multiplicity equal to two. We use the same\footnote{We think there is a typo somewhere in the paper by \cite{Zhang2014} and a minus is missing in front of the feedback gain. We added it to the presentation of $\mathbf{K}$ in \eqref{eq:k} in order to reproduce the results of the original paper with fidelity. Without that minus sign the system is unstable for any value of the delay.} feedback gain reported by \cite{Zhang2014}:
\begin{equation}
\mathbf{K}=-\left[\begin{array}{rrr}-0.2694&0.0402&-0.0899\\0.0386&-0.2857&-0.1238\end{array}\right].
\label{eq:k}
\end{equation}

With all this elements set up, we perform the stability analysis of two individual systems, one for each eigenvalue different than zero, following the discussions in section \ref{sec:ctcr}. The kernels of both subsystems have a size of three. For the system corresponding to $\lambda=1.3820$, the kernel delays are $1.3213$, $3.2785$ and $6.5033$ s. The root tendencies of the first two crossings are $+1$, indicating a destabilizing effect, whereas the last delay is stabilizing, it has $RT=-1$. The system corresponding to $\lambda=3.6180$ has kernel delays of value $0.9010$, $1.3971$ and $10.0999$ s, with root tendencies of $+1$, $+1$, and $-1$, respectively. Just from this results alone, we can conclude that the exact stability bound of the system is $\tau^{*}=0.901$ seconds. This is 2.5 times the value of 0.35 seconds reported by \cite{Zhang2014}. 

\begin{figure}
\centering
\includegraphics[scale=1]{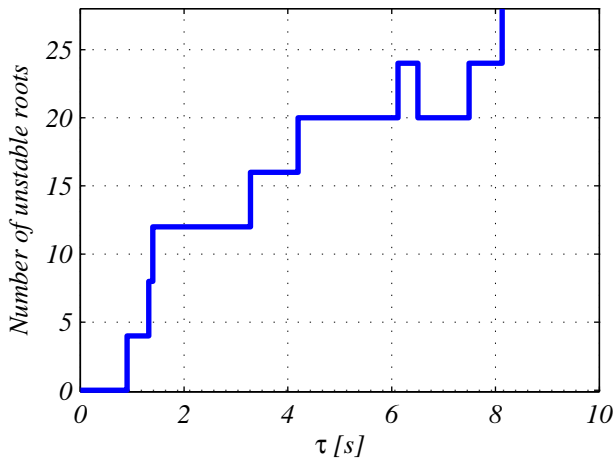}
\caption{Number of unstable roots of the system as a function of the delay.}
\label{fig:nu}
\end{figure}

To finalize the result, we deploy th D-subdivision concept and obtain $NU$, the number of unstable characteristic roots of the complete system as a function of the delay. This is shown in Figure \ref{fig:nu}. As we can see in this plot, there is only one interval for which $NU=0$. This is the only stable operating region for this case. At each transition point the number of unstable root changes by four because there are two subsystems created by each eigenvalue. Each subsystem contributes with two new roots at the transition point.

\begin{figure}
\centering
\subfigure[Stable behavior, $\tau=0.7$ s.]{\includegraphics[scale=0.9]{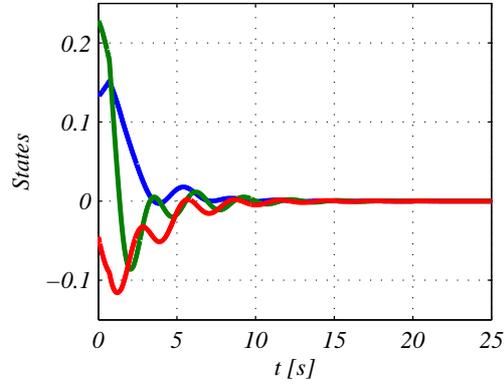}}
\subfigure[Marginally stable behavior, $\tau=0.9010$ s.]{\includegraphics[scale=0.9]{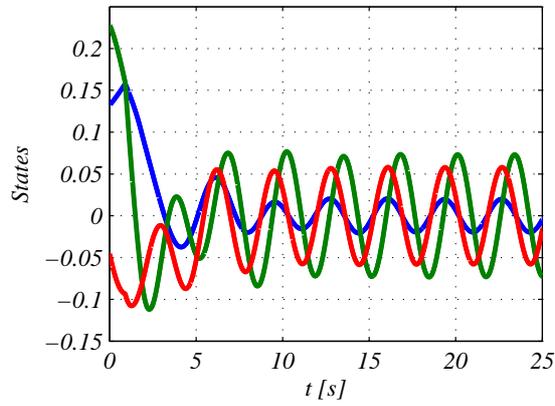}}
\subfigure[Unstable behavior, $\tau=1.1$ s.]{\includegraphics[scale=0.9]{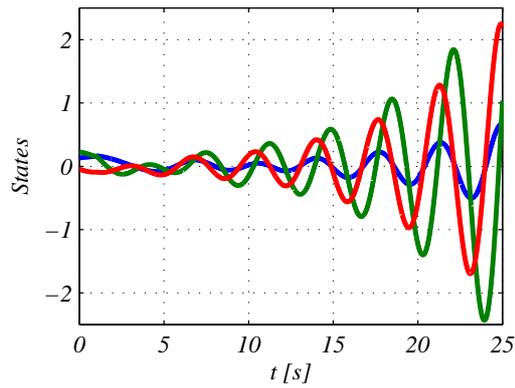}}
\caption{Simulation results for the disagreement subsystem corresponding to $\lambda=3.6180$ for different values of the delay.}
\label{fig:plots}
\end{figure}

To confirm the exactitude of the delay bound obtained we present simulation results with delay values smaller than the bound ($\tau=0.7$ seconds), exactly at the bound ($\tau=0.9010$ seconds) and larger than the bound ($\tau=1.1$ seconds). Rather than plotting the complete 15 states of the agents, we only show the results of the three states of the transformed disagreement subsystem \eqref{eq:subsys} corresponding to $\lambda=3.6180$, which is the one that introduces the instability first. This is done to avoid overcrowding the plot and to allow a better visualization of the qualitative features of the results, which would be obscured if over the disagreement states we plotted the average value, which is always increasing. The three plots can be seen in Figure \ref{fig:plots}

The code used in the examples is available at \texttt{https://db.tt/Ew1JWnES} or can be requested via email to the author.

\section{Discusion on Switching Topologies and the Most Exigent Eigenvalue}\label{sec:switch}
The work of \cite{Zhang2014} includes the analysis of the delay bound for a system in which the communication structure is switching among a set of $N$ known different topologies. Finding a delay bound for this case is as simple as it was for the fixed topology case. The eigenvalues of the $N$ Laplacian matrices are found, and the analysis of section \ref{sec:ctcr} is repeated for the different eigenvalues, which can be at most $N(n-1)$ once we rule out the zero eigenvalue common to all Laplacians. Since this analysis does not take much time once it is programmed, this is done very efficiently using CTCR, and the exact stability bound for the system under switching topologies can be found. In fact, for the example case on switching topologies presented by \cite{Zhang2014}, in which all the communication networks are star topologies with 50 agents, only two different eigenvalues must be analysed.

A further simplification of the task, however, can be achieved by studying the concept of \emph{most exigent eigenvalue}, introduced by \cite{IJC1}, for this class of systems. The most exigent eigenvalue is an eigenvalue of a certain matrix related to the communication topology that creates the most restrictive stability bound in the domain of the delay. For the consensus protocol studied by \cite{IJC1}, it is proven that the most exigent eigenvalue is the smallest eigenvalue of a weighted adjacency matrix of the system. For delayed consensus protocols described by Laplacians, \cite{TAC} observed that the most exigent eigenvalue is the largest eigenvalue of the Laplacian matrix. As we can see in Section \ref{sec:simulations}, this is also the case in the present paper, the largest eigenvalue of the Laplacian is also the most exigent. If this fact is proven strictly, a fantastic reduction in the complexity of the problem would be achieved: it would be sufficient to study only one system like \eqref{eq:subsys} to obtain an exact stability bound. This topic is a matter of further research. 

\section{Conclusions}\label{sec:conclusions}
This paper treated the problem of finding an exact bound for the delay that can be tolerated by a group of agents before it loses the ability to reach consensus. The agents are assumed to be driven by generic high order linear dynamics. We revisit a consensus algorithm recently proposed for this class of systems and revisit it using the Cluster Treatment of Characteristic Roots paradigm for its stability analysis.

The analysis consists in two steps. The firs one is a diagonalization of the systems, which decomposes it into individual systems of reduced order easier to analyse. The exact delay bound for each system is found using CTCR. The crucial first step of CTCR, the determination of the kernel and offspring sets, an approach based on the Kronecker summation operation is employed. The stability regions of all subsystems are intersected to obtain the allowable delay bound for the whole system.

When compared to other approaches to the same problem, the results of this work present the advantages of being exact (there is no conservatism), exhaustive (no stability region is undetected) and straightforward.

Avenues for further research include a detailed study of the most exigent eigenvalue for this class of systems and a more detailed exploration of the switching topologies case.

\bibliographystyle{plain}
\bibliography{swarms}

\begin{thebibliography}{10}

\bibitem{Bernstein2005}
D.S. Bernstein.
\newblock {\em Matrix Mathematics}.
\newblock Princeton University Press, Princeton, NJ., 2005.

\bibitem{Brewer1978}
J.W. Brewer.
\newblock Kronecker products and matrix calculus in system theory.
\newblock {\em IEEE Transactions on Circuits and Systems}, 25:772--781, 1978.

\bibitem{Cao2013b}
Y.~Cao, W.~Yu, W.~Ren, and G.~Chen.
\newblock An overview of recent progress in the study of distributed
  multi-agent coordination.
\newblock {\em IEEE Transactions on Industrial Informatics}, 9(1):427--438,
  2013.

\bibitem{IJC2}
R.~Cepeda-Gomez and N.~Olgac.
\newblock Consensus analysis with large and multiple communication delays using
  spectral delay space concept.
\newblock {\em International Journal of Control}, 84(12):1996--2007, 2011.

\bibitem{TAC}
R.~Cepeda-Gomez and N.~Olgac.
\newblock An exact method for the stability analysis of linear consensus
  protocols with time delay.
\newblock {\em IEEE Transactions on Automatic Control}, 56(7):1734--1740, 2011.

\bibitem{IJC1}
R.~Cepeda-Gomez and N.~Olgac.
\newblock Exhaustive stability analysis in a consensus system with time delay
  and irregular topologies.
\newblock {\em International Journal of Control}, 84(4):746--757, April 2011.

\bibitem{EJC}
R.~Cepeda-Gomez and N.~Olgac.
\newblock Stability analysis for the group dynamics consensus with time delay.
\newblock {\em European Journal of Control}, 18(5):456--468, September 2012.
\newblock 2012.

\bibitem{SCL}
R.~Cepeda-Gomez and N.~Olgac.
\newblock Exact stability analysis of second-order leaderless and
  leader-follower consensus protocols with multiple time delays.
\newblock {\em Systems and Control Letters}, 62(6):482--495, June 2013.

\bibitem{Ergenc2007}
A.~F. Ergenc, N.~Olgac, and H.~Fazelinia.
\newblock Extended kronecker summation for cluster treatment of lti systems
  with multiple delays.
\newblock {\em SIAM Journal on Control and Optimization}, 46(1):143--155, 2007.

\bibitem{Fazelinia2007}
H.~Fazelinia, R.~Sipahi, and N.~Olgac.
\newblock Stability robustness analysis of multiple time-delayed systems using
  building block concept.
\newblock {\em IEEE Transactions on Automatic Control}, 52(5):799--810, 2007.

\bibitem{Godsil}
C.~Godsil and G.~Royle.
\newblock {\em Algebraic Graph Theory}.
\newblock Springer, New York, 2004.

\bibitem{kolmanovski1986}
V.~B. Kolmanovski and V.~R. Nosov.
\newblock {\em Stability of Functional Differential Equations}.
\newblock Academic, London, U.K., 1986.

\bibitem{Lin2009-2}
J.~Lin and Y.~Jia.
\newblock Consensus of second-order discrete multi-agent systems with
  nonuniform time delays and dynamically changing topologies.
\newblock {\em Automatica}, 45(9):2154--2158, 2009.

\bibitem{Lin2009-1}
P.~Lin and Y.~Jia.
\newblock Further results on decentralised coordination in networks of agents
  with second-order dynamics.
\newblock {\em IET Control Theory and Applications}, 3(7):957--970, 2009.

\bibitem{Lin2008-1}
P.~Lin, Y.~Jia, J.~Du, and S.~Yuan.
\newblock Distributed control of multi-agent systems with second-order agent
  dynamics and delay-dependent communications.
\newblock {\em Asian Journal of Control}, 10(2):254--259, 2008.

\bibitem{Lin2008-2}
P.~Lin, Y.~Jia, and L.~Li.
\newblock Distributed robust $h_\infty$ consensus control in directed networks
  of agents with time-delay.
\newblock {\em Systems and Control Letters}, 57(8):643--653, 2008.

\bibitem{Liu2009}
X.~Liu and Y.-P. Tian.
\newblock Formation control of multi-agent systems with heterogeneous
  communication delays.
\newblock {\em International Journal of Systems Science}, 40(6):627--636, 2009.

\bibitem{Meng2011}
Z.~Meng, W.~Ren, Y.~Cao, and Y.~Zheng.
\newblock Leaderless and leader-follower consensus with communications and
  input delays under a directed network topology.
\newblock {\em IEEE Transactions on Systems, Man and Cybernetics - Part B},
  41(1):75--88, 2011.

\bibitem{Muenz2010}
U.~M\"unz, A.~Papachristodoulou, and F.~Allgower.
\newblock Delay robustness in consensus problems.
\newblock {\em Automatica}, 46(8):1252--1265, 2010.

\bibitem{Olfati-Saber2004}
R.~Olfati-Saber and R.M. Murray.
\newblock Consensus problems in networks of agents with switching topology and
  time-delays.
\newblock {\em IEEE Transactions on Automatic Control}, 49(9):1520--1533, 2004.

\bibitem{Olgac1994}
N.~Olgac and B.~Holm-Hansen.
\newblock A novel active vibration absorption technique: the delayed resonator.
\newblock {\em Journal of Sound and Vibration}, 176(1):93--104, 1994.

\bibitem{Olgac2002}
N.~Olgac and R.~Sipahi.
\newblock An exact method for the stability analysis of time-delayed linear
  time invariant systems.
\newblock {\em IEEE Transactions on Automatic Control}, 47(5):793--797, 2002.

\bibitem{Olgac2005}
N.~Olgac and R.~Sipahi.
\newblock Complete stability robustness of third order lti multiple time delay
  systems.
\newblock {\em Automatica}, 41(8):1413--1422, 2005.

\bibitem{Olgac2006}
N.~Olgac and R.~Sipahi.
\newblock An improved procedure in detecting the stability robustness of
  systems with uncertain delay.
\newblock {\em IEEE Transactions on Automatic Control}, 51(7):1164--1165, 2006.

\bibitem{Peng2009}
K.~Peng and Y.~Yang.
\newblock Leader-following consensus problems with a varying velocity leader
  and time-varying delays.
\newblock {\em Physica A: Statistical Mechanics and its applications},
  388(2-3):193--208, 2009.

\bibitem{Ren2008}
W.~Ren.
\newblock On consensus algorithms for double integrator dynamics.
\newblock {\em IEEE Transactions on Automatic Control}, 53(6):1503--1509, 2008.

\bibitem{Ren2005}
W.~Ren and R.~Beard.
\newblock Consensus seeking in multi-agent systems under dynamically changing
  interaction topologies.
\newblock {\em IEEE Transactions on Automatic Control}, 50(5):655--661, 2005.

\bibitem{kron}
R.~D. Schaefer.
\newblock {\em An Introduction to Nonassociative Algebras}.
\newblock Dover, 1996.

\bibitem{Sun2009-2}
Y.~Sun and L.~Wang.
\newblock Consensus of multi-agent systems in directed networks with nonuniform
  time-varying delay.
\newblock {\em IEEE Transactions on Automatic Control}, 54(7):1607--1613, 2009.

\bibitem{Sun2009-1}
Y.~Sun and L.~Wang.
\newblock Consensus problems in networks of agents with double integrator
  dynamics and time delays.
\newblock {\em International Journal of Control}, 82(10):1937--1945, 2009.

\bibitem{Vicsek1995}
T.~Vicsek, A.~Czir\'ok, E.~Ben-Jacob, I~Cohen, and O.~Shochet.
\newblock Novel type of phase transition in a system of self-driven particles.
\newblock {\em Physical Review Letters}, 75(6):1226--1229, 1995.

\bibitem{Zhang2014}
Ya~Zhang and Yu-Pin Tian.
\newblock Allowable delay bound for consensus of linear multiagent systems with
  communication delay.
\newblock {\em International Journal of Systems Science}, 45(10):2172--2181,
  2014.

\end{thebibliography}

\section*{Appendix A: The Kronecker Summation}\label{sec:appendix}
Consider two square matrices $\mathbf{A}\in\mathbb{R}^{n\times n}$ and $\mathbf{B}\in\mathbb{R}^{m\times m}$. The Kronecker summation of these two matrices, denoted by $\oplus$, is defined as:
\begin{equation}
\mathbf{M}=\mathbf{A}\oplus\mathbf{B}=\mathbf{A}\otimes\mathbf{I}_m+\mathbf{I}_n\otimes\mathbf{B},
\label{eq:kronsum}
\end{equation}
where $\otimes$ denotes the Kronecker product. The result or this operation is a square matrix of dimensions $nm\times nm$. Some extra discussions on this operation are found in the works of \cite{Brewer1978} and \cite{Bernstein2005}. 

For the purpose of this work, the most interesting property of the Kronecker summation is that the $nm$ eigenvalues of $\mathbf{M}$ are indeed pairwise additions of the eigenvalues of $\mathbf{A}$ and $\mathbf{B}$. That is, if the eigenvalues of $\mathbf{A}$ are $\lambda_{A1},\,\lambda_{A2},\,\ldots,\,\lambda_{an}$ and those of $\mathbf{B}$ are $\mathbf{A}$ are $\lambda_{B1},\,\lambda_{B2},\,\ldots,\,\lambda_{Bm}$, the eigenvalues of $\mathbf{M}$ are $\lambda_{M1}=\lambda_{A1}+\lambda_{B1},\ \lambda_{M2}=\lambda_{A1}+\lambda_{B2},\,\ldots\,\lambda_{M,{mn}}=\lambda_{An}+\lambda_{Bm}$. The Kronecker summation, in fact, is also known as \emph{eigenvalue addition}.
\end{document}